\documentclass[aoas,preprint,12pt]{article}
\RequirePackage[OT1]{fontenc}
\RequirePackage{amsthm,amsmath}
\usepackage{soul}
\usepackage[numbers,sort&compress]{natbib}
\RequirePackage[colorlinks=black,citecolor=blue,urlcolor=blue]{hyperref}
\usepackage{amsfonts}
\usepackage{amssymb}
\usepackage{graphicx}
\usepackage{setspace}
\newtheorem{theorem}{Theorem}[section]
\newtheorem{lemma}{Lemma}[section]

\newtheorem{example}{Example}[section]
\newtheorem{definition}{Definition}[section]
\newtheorem{problem}{Problem}[section]
\numberwithin{equation}{section}
\newcommand{\bxi}{{\boldsymbol \xi}}
\newcommand{\bkappa}{{\boldsymbol \kappa}}

\newcommand{\bet}{{\boldsymbol \eta}}
\newcommand{\bom}{{\boldsymbol \omega}}

\newcommand{\blambda}{{\boldsymbol \lambda}}
\newcommand{\bzeta}{{\boldsymbol \zeta}}
\def\by{{\boldsymbol y}}

\def\b1{{\mathbf 1}}
\def\bw{{\boldsymbol w}}
\def\bA{{\boldsymbol A}}

\def\bC{{\boldsymbol C}}
\def\bu{{\boldsymbol u}}
\def\bz{{\boldsymbol z}}
\def\bh{{\boldsymbol h}}
\def\be{{\boldsymbol e}}
\def\bk{{\boldsymbol k}}
\def\bv{{\boldsymbol v}}

\def\bp{{\boldsymbol p}}
\def\bq{{\boldsymbol q}}

\def\bx{{\boldsymbol x}}

\def\bI{{\boldsymbol I}}

\def\bH{{\boldsymbol H}}

\def\bW{{\boldsymbol W}}

\def\cM{{\mathcal M}}
\def\cK{{\mathcal K}}
\def\cF{{\mathcal F}}
\def\cP{{\mathcal P}}

\def\bbR{{\mathbb R}}

\def\bbP{{\mathbb P}}
\def\bbQ{{\mathbb Q}}

\def\btau{{\boldsymbol \tau}}
\def\blambda{{\boldsymbol \lambda}}

\graphicspath{%
    {converted_graphics/}
    {/}
}

\begin{document}

\title{Joint probabilities under expected value constraints, transportation problems, maximum entropy in the mean, and geometry in the space of probabilities} 
\author{Henryk Gzyl\\
Centro de Finanzas IESA, Caracas, Venezuela.\\
 henryk.gzyl@iesa.edu.ve}

\date{}
 \maketitle

\setlength{\textwidth}{4in}

\vskip 1 truecm
\baselineskip=1.5 \baselineskip \setlength{\textwidth}{6in}
\begin{abstract}
There are interesting extensions of the problem of determining a joint probability with known marginals. On the one hand, one may impose size constraints on the joint probabilities. On the other, one may impose additional constraints like the expected values of known random variables.

If we think of the marginal probabilities as demands or supplies, and of the joint probability as the fraction of the supplies to be shipped from the production sites to the demand sites, instead of joint probabilities we can think of transportation policies. Clearly, fixing the cost of a transportation policy is equivalent to an integral constraints upon the joint probability.

We will show how to solve the cost constrained transportation problem by means of the method of maximum entropy in the mean. We shall also show how this approach leads to an interior point like method to solve the associated linear programming problem. We shall also investigate some geometric structure the space of transportation policies, or joint probabilities or pixel space, using a Riemannian structure associated with the dual of the entropy used to determine bounds between probabilities or between transportation policies.

\end{abstract}


\begin{spacing}{0.5}
\small{\tableofcontents}
\end{spacing}

\section{Introduction and Preliminaries} 
In a previous paper Gzyl (2020) developed an approach to the problem of determining a joint probability on $\{1,...,N\}\times\{1,...,N\}$ when, besides the specification of its marginals, we might have constraints on the range of values of the joint probability in a given cell. Besides these constraints, we might have additional information specified as the expected value of a collection of random variables with respect to the unknown probability.

An interesting twist on the problem of comparing two histograms on the set $\{1,...,N\}$ is to think of them as mounds of dirt to be transformed one onto another, and amount of dirt taken from one of the mounds to the other, is described by a function on $\{1,...,N\}^2$ whose marginals are the amount of dirt (the given histograms) at each point.  Therefore, a procedure to transform one dirt distribution onto the other at a given cost is of clear interest. And the smaller the cost, the better.

 Two comprehensive textbooks about the transportation problem in the continuous case plus a guide to a large body of literature are Santambrogio's (2015) and Villani's (2008). We direct the reader to Chapter 13 of Kapur's (1998) for early applications of the standard method of maximum entropy to the discrete transportation problem.

To introduce notations for the discrete case, let $\{p_i: i=1,...,N\}$ and $\{q_j: j=1,...,N\}$ on  $\{1,...,N\}$ be two probability assignments on $\{1,...,N\}.$ There are infinitely many joint probabilities $\pi_{i,j}:, i,j=1,...,N$ on  $\{1,...,N\}^2$ such that 
\begin{equation}\label{const1}
\sum_{j=1}^N\pi_{i,j} = p_i,\;\;\;\sum_{i=1}^N\pi_{i,j} = q_j.
\end{equation}
Therefore, besides designing a method to solve the problem, we need criteria to  choose among solutions. 
It is interesting, specially for numerical purposes, to relabel the unknowns and the constrains. For that we list the elements of  $\{1,...,N\}\times\{1,...,N\}$ in lexicographic order. To he pair $(i,j)$ we associate $n=(i-1)N+j$ and to the joint probability $\pi$ we associate a $N^2-$vector $\bx\in [0,1]^{N^2}.$ In the first appendix we explain how to rewrite the constraints in matrix form as:   
\begin{equation}\label{const1.1}
\bC\bx = {\bp \atopwithdelims[]\bq} = \by.
\end{equation}
\noindent  where we put $\bp=(p_1,...,p_N)^t$ and $\bq=(q_1,...,q_N)^t,$ where the superscript $t$ denotes transposition (we think of vectors as columns). We shall also see in the first appendix that $\bC$ is not of full rank. 

\subsection{Extension of the problem of reconstruction of a probability from its marginals}
To state the transportation problem, we think of $\pi_{i,j}$ as a fraction of ``goods'' being ``transported'' from point $i$ to point $j,$ and denote by $W_{i,j}$ the cost of doing that. The cost of transforming the probability vector $\bp$ onto the probability vector $\bq$ is modeled by a random variable $\bW$ is $\sum_{i,j}W_{i,j}\pi_{i,j}$ which after relabeling becomes $\sum_nW_n\pi_n$ which can be thought of as the expected value of the random variable $\bW.$ When a cost to be met is fixed, to the constraints (\ref{const1.1} we must add the cost constraint which we write as
\begin{equation}\label{cconst}
\sum_{n=1}^{N^2} W_{n}\pi_{n} = c.
\end{equation}
To consider (\ref{const1.1}) and (\ref{cconst}) in a unified way augment the constraint matrix $\bC$ by adding $\bW^t$ as a last row to it, and denote the new matrix by $\bA.$ The data vector is extended by adding the cost constraint $c$ as its $2N+1-$th component.

After relabeling, the problem to solve to recover the joint probability $\pi_{i,j}$ becomes:
\begin{problem}\label{prob2}
Determine $\bx(c)^*\in[0,1]^{N^2}$ such that
\begin{equation}\label{prob2.1}
\bA\bx(c)^* = \by(c)
\end{equation}
where explicitly 
\begin{equation}\label{augm1}
 \bA ={\bC \atopwithdelims[]\bW^t},\;\;\;\by(c) = {\by \atopwithdelims[] c}.
\end{equation}
\end{problem} 

We use a generic $c$ as argument in the augmented vector because we are going to be varying that $c$ as part of the application of maximum entropy in the mean to the transportation problem. To close this section, we mention that the standard transportation problem can be stated as:

\begin{problem}\label{prob1}
Determine $\bx^*\in[0,1]^{N^2}$ such that
\begin{equation}\label{prob1.1}
\bx^* = argmin\{\langle\bW,\bx\rangle: \bC\bx=\by\}
\end{equation}
where $\bC$ and $\by$ were introduced in (\ref{const1}) after the relabeling. 
\end{problem}

\subsection{Size constraints upon the solution}\label{sizecon}
It may happen that one has prior information upon the solution in the form of size constraints. In this case Problem \ref{prob2} is to be replaced by
\begin{problem}\label{prob2a}
Determine $\bx(c)^*\in\prod_{n=1}^{N^2}[a_n,b_n]$ such that
\begin{equation}\label{prob2a.1}
\bA\bx(c)^* = \by(c)
\end{equation}
\end{problem}
\noindent where the rest of the symbols are as in (\ref{augm1}). An interpretation of the constraints from the demand-consumption point of view could be the following. If $n\Leftrightarrow (i,j),$ then $x_n\in[a_n,b_n],$ may mean that site $i$ requires at least $a_{i,j}$ and at most $b_{i,j}$ units of some good from site $j.$

\subsection{Comments about the unconstrained inverse problem}
The solution to the linear algebraic equation
$$\bA\bx(c)^* = \by(c)$$
\noindent is given by 
$$\bx = \bA^+\by(c) + \big(\bI -\bA^+\bA\big)\bz,$$
\noindent where $\bz$ is any arbitrary element in $\bbR^{N^2}$ and $\bA^+$ denotes the Moore-Penrose inverse of $\bA.$  Note that $\bA^+\bA$ is a projection onto $Ker(\bA)^{\perp},$ and therefore $\big(\bI -\bA^+\bA\big)\bz$ is in $Ker(\bA).$\\
The whole difficulty in Problem \ref{prob2a} lies in how to choose one of the infinitely many solutions that also satisfies the convex constraints $x_n\in[a_n,b_n].$ To conclude this section we add that the method maximum entropy in the mean (MEM) is specially suited to deal 
with this.

\subsection{Undoing the relabeling}
Let us now denote by $x^*$ the solution that we are after. To write is as $x^*_{i,j}$ we proceed as follows. Write $n=kN+r$ with $r=0,1,...,N-1.$ Then
\begin{eqnarray}\label{undo}\nonumber
\mbox{If}\;\;r=0,\;\;\mbox{then}\;\; (i,j)=(k+1,N).\\\nonumber
\mbox{If}\;\;1\leq r\leq N-1,\;\;\mbox{then}\;\; (i,j)=(k+1,r)
\end{eqnarray}

\subsection{Organization of the paper}\label{contents}
Section 2 is devoted to the fixed cost transportation problem described in (\ref{prob2a.1}). To begin with, in Section 2 we explain the basics of MEM and how it is applied to solve that problem. In Section 3 we provide a way to compare probability laws defined on $\{1,...,N\}\times\{1,...,N\}$ or equivalently, on $\{1,...,N^2\}$ constrained to take values in $\prod(a_n,b_n).$ For that we shall define a Riemannian metric on  $\prod(a_n,b_n)$ by pulling back a metric on $\bbR^{N^2}$ obtained as the Hessian of $\bzeta(\btau)=\sum\ln\big(e^{-a_n\tau_n}+e^{-b_n\tau_n}\big).$ Observe that for $\btau=\bA^*\blambda$ we obtain $Z(\blambda)=\bzeta(\bA^*\blambda).$ This is why we start form $\bzeta(\btau)$ to define the Riemann metric in $\bbR^{N^2}.$

In Section 4 we explain how the maxentropic approach could be applied to solve the minimum cost transportation problem.  The procedure consists of decreasing the cost $c$ until the corresponding joint probability is a point at the boundary of the probability simplex, we will see that when the data vector approaches the boundary of $\bA(\cP),$ the solution to the maxentropic problem ceases to exist. In particular, we shall understand why the numerical procedure breaks down when that happens. The presentation is based on Gamboa and Gzyl (1990) which was first attempt to a maxentropic approach to the interior point approach to finite dimensional linear programming.

\section{The MEM approach to solve Problems (\ref{prob2.1}-\ref{prob2a.1})}\label{MEM}
The standard method of maximum entropy proposed by Jaynes (1957) to solve a problem in statistical physics, consisting of determining a density from the knowledge of the expected value of a few random variables. Such densities characterize thermal equilibrium in statistical thermodynamics. About the same time Kullback proposed the method to solve a similar problem in statistics.  Campbell (1966) interpreted the method from the point of view of maximum likelihood, and Good (1963) used the method of maximum entropy to determine the joint probability of a pair of discrete valued random variables when only their marginal probability are known. The aim of Gzyl's (2020) was to use maximum entropy in the mean (MEM) to take care of natural constraints (prior information on the range) on the unknown probabilities.

Here we sketch essence the method here and direct the reader to Dacunha-Castelle and Gamboa (1990). Our presentation is along the lines developed in Golan and Gzyl (2002).

To describe MEM and the solution to (\ref{prob2a}) we need to introduce some notation. (Consider an auxiliary probability space $(\Omega,\cF,\bbP)$ where:
$$\Omega = \prod_{n=1}^{N^2}[a_n,b_n]$$
and $\cF$ denotes the Borel subsets of $\Omega,$ and $\bbP$ is a (yet unspecified) probability measure on $(\Omega,\cF).$ Denote by $\bom$ the points in $\Omega$ and let $\bxi:\Omega\to\Omega$ denote the identity mapping. Note that $\bxi_n(\bom) \in [a_n,b_n]\;\;\mbox{for}\;\;1\leq n\leq N^2.$ 

 With all this, the MEM procedure consists of replacing the constrained linear problem (\ref{prob2}) by the following problem:
\begin{problem}\label{prob4}
Determine a probability $\bbP$ measure on $(\Omega,\cF)$ such that
\begin{equation}\label{prob4.1}
\bA E_\bbP[\bxi] = \by(c)
\end{equation}
\noindent where, recall, $\bA$ and $\by(c)$ are given by
\begin{equation}\label{augm1.a}
 \bA ={\bC \atopwithdelims[]\bW^t},\;\;\;\by(c) = {\by \atopwithdelims[] c}.
\end{equation}
\end{problem}

{\bf Comment:} Notice that if such a $\bbP$ is found, then $\bx^*=E_\bbP[\bxi]$ satisfies problem (\ref{prob4}). Notice that rendering problem (\ref{prob2a.1}) as problem \ref{prob4.1} automatically takes care of the constrains. This is the essence of MEM.

In Section \ref{MD} we explain how to obtain $\bbP.$ In the setup that we use it turns out that
\begin{equation}\label{sol1g}
x^*_n =  \frac{a_ne^{-a_n(\bA^t\blambda^*)_n }+b_ne^{-b_n(\bA^t\blambda^*)_n }}{e^{-a_n(\bA^t\blambda^*)_n }+e^{-b_n(\bA^t\blambda^*)_n }}.
\end{equation}
This was the generic case treated in Gzyl (2020). When $a_n=0, b_n=1$ for all $n,$ the notations simplify an the solution to (\ref{prob2.1}) is given by
\begin{equation}\label{sol1}
x^*_n =  \frac{e^{-(\bA^t\blambda^*)_n }}{1+e^{-(\bA^t\blambda^*)_n }}.
\end{equation}
\noindent where $\blambda\in \bbR^{2N+1}$ is a vector of Lagrange multipliers.
It couldn't be more clear: the components of $\bx$ are convex combinations of the end point values that define the constraints, thus the constraints are met. In our case, the points are $\{a_n,b_n\}$ for $1\leq n\leq N^2.$ It only remains to mention that $\blambda^*$ is determined minimizing the strictly convex function defined on $\bbR^{2N+1}$ by
\begin{equation}\label{dualent1}
\Sigma(\blambda,\by(c)) = \ln Z(\blambda) + \langle\blambda,\by(c)\rangle.
\end{equation}
In our setup the function $Z(\blambda)$ is computed in Section \ref{MD} to be
\begin{equation}\label{norm1}
Z(\blambda) =\prod_{n=1}^{N^2}\big(e^{-a_n(\bA^t\blambda^*)_n } + e^{-b_n(\bA^t\blambda)_n}\big).
\end{equation}
Actually, using (\ref{norm1}) the representations (\ref{sol1}) can be read off form the first order condition for $\blambda^*$ to be a minimizer of (\ref{dualent1}).

\section{Bounds on the distances between probabilities}
In this section we are going to define a Riemannian metric on $\prod(a_n,b_n)$ by pulling back a Riemannian metric defined on $\bbR^{N^2}.$ We follow this route because the Riemannian metric on $\bbR^{N^2}$ is related to the Laplace transform of the reference measure $\bbQ$ that enters in the definition of the dual entropy function.

\subsection{Distance in $\bbR^{N^2}$}
Our starting point is the Laplace transform of the reference measure $\bbQ$ defined by
\begin{equation}\label{LT}
\zeta(\btau) = E_{\bbQ}[e^{\langle\btau,\bxi\rangle}] = \prod_{n=1}^{N^2}\bigg(e^{-a_n\tau_n}+e^{-b_n\tau_n}\bigg) = \prod_{n=1}^{N^2}\zeta_n(\tau_n),
\end{equation}
and then the moment generating function
\begin{equation}\label{LT1}
M(\btau) = \ln\zeta(\btau) = \sum_{n=1}^{N^2}\ln\zeta_n(\tau_n) = \sum_{n=1}^{N^2}\ln\bigg(e^{-a_n\tau_n}+e^{-b_n\tau_n}\bigg).
\end{equation}

To continue, with the notations introduced above, the Hessian matrix of the moment generating function $M(\btau)$ is diagonal with entries 
$$M_n''(\tau_n)=\frac{a^2_ne^{-a_n\tau_n}+b^2_ne^{-b_n\tau_n}}{e^{-a_n\tau_n}+e^{-b_n\tau_n}}-\bigg(\frac{a_ne^{-a_n\tau_n}+b_ne^{-b_n\tau_n}}{e^{-a_n\tau_n}+e^{-b_n\tau_n}}\bigg)^2.$$
Here $M'(\tau_n), M_n''$ denote the first and second derivative of $M_n$ with respect of its argument, and below we use the conventional $\dot{\tau}, \ddot{\tau}$ to denote first and second derivatives with respect to time . These are the diagonal elements of the covariance matrix of $\bxi$ with respect to the maxentropic probability obtained in Section \ref{MD}. After a simple calculation we can write
\begin{equation}\label{step1}
M_n''(\tau_n)=\bigg(\frac{b_n-a_n}{e^{(b_n-a_n)\tau_n/2}+e^{-(b_n-a_n)\tau_n/2}}\bigg)^2 = \bigg(\frac{(b_n-a_n)e^{-(b_n-a_n)\tau_n/2}}{e^{(b_n-a_n)\tau_n/2}+e^{-(b_n-a_n)\tau_n/2}}\bigg)^2 =\bigg(h_n'(\tau_n)\bigg)^2.
\end{equation}
\noindent where clearly 
$$h'_n(\tau_n) = 2\frac{d}{d\tau_n}\arctan(e^{(b_n-a_n)\tau_n/2}).$$
Since $h_n(\tau_n)=2\arctan\big[\exp\big((b_n-a_n)\tau_n/2\big)\big]$ is strictly increasing and continuously differentiable we can define the change of variables in $\bbR^{N^2}$ by $\bh(\btau)_n = h_n(\tau_n).$  With all this, the (square) of the velocity along a curve $t\to\btau(t)$ in the Hessian metric is given by
\begin{equation}\label{scapr}
\langle\dot{\btau},M''(\btau)\dot{\btau}\rangle = \langle\frac{d}{dt}\bh(\btau(t)),\frac{d}{dt}\bh(\btau(t))\rangle = \langle\dot{\bw},\dot{\bw}\rangle.
\end{equation}
Here we used the fact that $M_n''(\tau_n)=h'_n(\tau_n)h'_n(\tau_n)$ and we put $\bw(t)=\bh(\btau(t))$ for the curve in the $\bh$ coordinates. 
Therefore, in the $\bw$ coordinates the geodesics are straight lines given by $\bw(t)=\bw(0) + t\bkappa.$ To determine $\bw(0)$ and $\bkappa=\bw(1)-\bw(0)$ we suppose that the geodesic starts at $\btau(1)$ at $t=0$ and passes through $\btau(2)$ at $t=1.$ Therefore $\bw(0)=\bh(\btau(1)$ and $\bkappa=\bh(\btau(2)-\bh(\btau(1).$ Let $\bH$ denote the compositional inverse of $\bh$ with components $H_n=(h_n)^{-1}.$  With this $\btau(t)=\bH\big(\bh(\btau(1)+t\bkappa\big).$
All of this was established in Gzyl (2020b). We sum it up as:
\begin{theorem}\label{geod1}
In $\bbR^{N^2}$ consider the Riemannian distance defined by the Hessian of $M(\btau).$ Let $\btau(0)$ and $\btau(1)$ be any two points in $\bbR^{N^2}.$ The the geodesic that starts at $\btau(0)$ at $t=0$ and arrives at $\btau(1)$ at $t=1$ solves:
\begin{equation}\label{eqgeo1}
M_n''(\tau_n)\ddot{\tau}_n + M_n(\tau_n)\big(\dot{\tau}_n\big)^2, \;\;\;n=1,...,N^2.
\end{equation}
The solution to these is shown to be:
\begin{equation}\label{geod2}
\tau_n(t) = H_n\bigg(h_n(\tau_n(0)) +\kappa_n t\bigg),\;\;\;\;n=1,....,N^2.
\end{equation}
\end{theorem}
\noindent where, as indicated above, $h_n(\tau_n)=2\arctan(e^{(b_n-a_n)\tau_n/2}),$ and we use $H_n=h_n^{-1}.$ Also, $\kappa_n=h_n(\tau_n(1)-h_n(\tau_n(0).$
It is also proved in Gzyl (2020b) that the geodesic distance between $\btau(0)$ and $\btau(1)$ is given by
\begin{equation}\label{geodis1}
d^2_M(\btau(0),\btau(1)) = \sum_{n=1}^{N^2}\big(h_n(\tau_n(1))-h_n(\tau_n(0))\big)^2.
\end{equation}

\subsection{The geometry on the pixel space}
Here we provide a geometry  in the space of parameters different from that considered by Amari et al. (2018). Recall that in order to obtain the solution to the maxentropic procedure we computed
\begin{equation}\label{sol1g.3}
\xi_n(\tau_n) = -\frac{d}{d\tau_n}M_n(\tau_n) = -\frac{d}{d\tau_n}\ln\big(e^{-a_n\tau_n}+e^{-b_n\tau_n}\big) = \frac{a_ne^{-a_n\tau_n }+b_ne^{-b_n\tau_n}}{e^{-a_n\tau_n}+e^{-b_n\tau_n}},
\end{equation}
\noindent and evaluated at $\tau_n = (\bA^t\blambda^*)_n,$ to obtain $x^*_n = \xi_n((\bA^t\blambda^*)_n).$ The passage from $\bbR^{N^2}$ to $\prod(a_n,b_n)$ is provided by
\begin{lemma}\label{corresp}
With the notations just introduced, the mapping defined by (\ref{sol1g.3}):
\begin{equation}\label{corresp1}
\bbR^{N^2}\to\prod_{n=1}^{N^2}(a_n,b_n),\;\;\;\btau \longleftrightarrow \bxi(\btau)
\end{equation}
\noindent is a continuously differentiable bijection.
\end{lemma}

With this mapping we can pull back the diagonal metric on $\bbR^{N^2}$ to a diagonal metric on $\prod(a_n,b_n)$ by setting 
\begin{equation}\label{PB0}
G_n(\xi_n) = M''(\tau_n(\xi_n)\big(\frac{d\tau_n}{d\xi_n}\big)^2.
\end{equation}
For not to repeat arguments similar to those in the previous section, we postpone the derivation of the equations of the geodesics the fact that they are the pullback of the geodesics of $M''(\btau),$ etc., to the last appendix. For the time being we relate the bounds on the distances between the geodesics in $\bbR^{N^2}$ to the bounds on the geodesics in $\prod(a_n,b_n)$ and the distance between solutions to the transportation problem. We proceed as follows. Start from (\ref{sol1g.3}) and invoke (\ref{step1}) along the way, note that, for any $\tau_n(1),\tau_n(2):$
\begin{eqnarray}
&\xi_n(\tau_n(1))-\xi_n(\tau_n(2)) = \frac{d}{d\tau_n}M_n(\tau_n(2))-\frac{d}{d\tau_n}M_n(\tau_n(1))= \int_{\tau_n(1)}^{\tau_n(2)}M''_n(u)du\nonumber\\
&= \int_{\tau_n(1)}^{\tau_n(2)}(h'(u))^2du = \int_{\tau_n(1)}^{\tau_n(2)}h'(u)dh(u) = \int_{\tau_n(1)}^{\tau_n(2)}\bigg(\frac{b_n-a_n}{e^{(b_n-a_n)\tau_n/2}+e^{-(b_n-a_n)\tau_n/2}}\bigg)dh(u).\nonumber
\end{eqnarray}
For the last step we replaced $h_n'(u)$ by its explicit form to note that, since $h_n'(u)\leq (b_n-a_n)/2,$ then
$$|\xi_n(\tau_n(1))-\xi_n(\tau_n(2))| \leq  \frac{(b_n-a_n)}{2}|h_n(\tau_n(1))-\tau_n(2)|.$$
From this a few inequalities are clear.
\begin{theorem}\label{compdis}
With the notations introduced above, let $L=\sup(b_n-a_n).$ Then
\begin{eqnarray}
\|\bxi(\btau(1))-\bxi(\btau(2)\|_2 \leq \frac{L}{2}d_M(\btau(1),\btau(2))\label{comps1}\\
\frac{1}{N}\|\bxi(\btau(1))-\bxi(\btau(2)\|_1 \leq  \frac{L}{2}d_M(\btau(1),\btau(2))\label{comps2}\\
\sup_n|\xi_n(\tau_n(1))-\xi_n(\tau_n(2))| \le \frac{L}{2}\sup_n|h_n(\tau_n(1))-h_n(\tau_n(2))|.\label{comps3}
\end{eqnarray}
\end{theorem}

In the last appendix we will show as well that
$$d^2(\bxi(2),\bxi(1)) = \sum_{n=1}^{N^2}\big(k_n(\xi_n(2)-k_n(\xi_n(1)\big)^2$$
where $\bk(\bxi)=\bh(\btau(\bxi))$ which results in $d_G(\bxi(1),\bxi(2))=d_M(\btau(1),\btau(2))$. This allows us to restate the last theorem as:

\begin{theorem}\label{compdisalt}
With the notations introduced for Theorem \ref{compdis}, let $d_G^2(\bxi(1),\bxi(2))=\sum_{n=1}^{N^2}\big(k_n(\xi_n(1))-k_n(\xi_n(2)\big)^2,$ Then
\begin{eqnarray}
\|\bxi(\btau(1))-\bxi(\btau(2)\|_2 \leq \frac{L}{2}d_G(\bxi(1),\bxi(2))\label{compsalt1}\\
\frac{1}{N}\|\bxi(\btau(1))-\bxi(\btau(2)\|_1 \leq  \frac{L}{2}d_G(\bxi(1),\bxi(2))\label{compsalt2}\\
\sup_n|\xi_n(\tau_n(1))-\xi_n(\tau_n(2))| \le \frac{L}{2}\sup_n|k_n(\xi_n(1))-k_n(\xi_n(2))|.\label{compsalt3}
\end{eqnarray}
\end{theorem}  

Observe that from $$\xi_n(\tau_n(1))-\xi_n(\tau_n(2)) = \int_{\tau_n(1)}^{\tau_n(2)}(h'(u))^2du$$
we obtain the bound $\|\bxi(\btau(1))-\bxi(\btau(2))\|_2 \leq \frac{L}{2}\|\btau(1)-\btau(2)\|_2$
which is a much worse bound than (\ref{comps1}) say, because there the distance on the right hand side is between points in the (bounded) range of $h.$ 

\subsection{Distance between solutions to the transportation problem}
Now we bring in the fact that the maxentropic solution to the transportation problem is obtained by setting $\bx(\blambda) = \bxi(\bA^*\blambda).$ The former bounds yield:

\begin{theorem}
Suppose that we have only one cost constraint, then the dimension of the space of constraints is $K=2N+1.$ Actually, when $\blambda$ sweeps $\bbR^{K},$ then $\bx^*(\blambda)=\bxi(\bA^*\blambda)$ ranges over the solutions to problem (\ref{prob2a}) with constraint $\bA\bxi(\bA^*\blambda).$ For $\blambda(1),\blambda(2)\in\bbR^{K}$ we have
\begin{eqnarray}
\|\bx^*(\blambda(1))-\bx^*(\blambda(2))\|_2 \leq \frac{L}{2}d_M(\bA^*\blambda(1),\bA^*\blambda(2))\label{comps1.1}\\
\frac{1}{N}\|\bx^*(\blambda(1))-\bx^*(\blambda(2))\|_1 \leq  \frac{L}{2}d_M(\bA^*\blambda(1),\bA^*\blambda(2))\label{comps2.1}\\
\sup_n |\bx^*(\blambda(1))-\bx^*(\blambda(2))| \le \frac{L}{2}\sup_n|h_n((\bA^*\blambda)_n(1))-h_n(\bA^*\lambda)_n(2)|.\label{comps3.1}
\end{eqnarray}
\end{theorem}

To conclude, we shall verify that the right hand side  of, say (\ref{comps1}), is the geodesic distance between the two points. For that, let us now write $\Phi(\blambda)=\ln Z(\blambda),$ that is,
$$\Phi(\blambda) = \ln Z(\lambda)=\ln\zeta(\bA^*\blambda)=\sum_{n=1}^{N^2}\ln\bigg(e^{-a_n(\bA^t\blambda^*)_n }+e^{-b_n(\bA^t\blambda^*)_n }\bigg).$$
The Hessian matrix of $\Phi$ in terms of that of $M(\tau)$ is $Hess\Phi(\blambda)=\bA HessM(\bA^*\blambda)\bA^*.$ That is, the metric is not diagonal anymore, but we can use the change of variables introduced above to verify the theorem.  Note that if $t \to \blambda(t)$ is a curve in $\bbR^{2N+1}$parameter space, using the change of variables $\bv=\bh(\bA^*\blambda),$ we can compute the length of the velocity $\dot{\blambda}(t)$ as
$$\langle\dot{\blambda},Hess\Phi(\blambda)\dot{\blambda}\rangle = \langle\bA^*\dot{\blambda},HessM(\bA^*\blambda)\bA^*\dot{\blambda}\rangle = \langle\dot{\bv},\dot{\bv}\rangle.$$ 
The curves that minimize the distance in the $\bv$ coordinates are straight lines. Since the solutions to the transportation are given by $\bx^*=\bxi(\bA^*\blambda),$ solving for $\bA^*\blambda$ using the change of variables is enough for our purposes. With end conditions $\blambda(1)$ at $t=0$ and $\blambda(2)$ at $t=1.$ From the comments made above, we have 
\begin{equation}\label{inter}
\bA^*\blambda(t) = \bH\bigg(\bh(\bA^*\blambda(1)) + t\bet\bigg)
\end{equation}
\noindent where here $\bet=\bh(\bA^*\blambda(2))-\bh(\bA^*\blambda(1)).$ It is therefore clear that we can interpret the right hand sides of (\ref{comps1})-(\ref{comps3}) as geodesic distances. We mention in passing, that the most that we can say is that, up to a term in $Ker(\bA^*),$ we have
\begin{equation}\label{sollam}
\blambda(t) = (\bA^*)^+\bH\bigg(\bh(\bA^*\blambda(1)) + t\bet\bigg)
\end{equation}
\noindent where $(\bA^*)^+$ is the Moore-Penrose inverse of $\bA^*.$ 

To finish, we mention the following addendum to Gzyl (2020). When there are no constraints besides the marginals, and if $a_n=0, b_n=1$ for all $n=1,...,N^2,$ then (\ref{comps1})-(\ref{comps3}) compares standard distances between two probability densities versus the geodesic distance between them obtained from the geometric induced by the moment generating function.

\section{The entropic approach to the minimal cost transportation problem}
As we said in Section (\ref{contents} the procedure is to apply MEM to a sequence of decreasing fixed cost problems. For that, we must first prescribe a method to choose an initial point and the how to decrease the step and when to stop.
\subsection{Choosing an initial cost $c_0$}
To fix an initial point for the iterative process, use the marginals to define the product probability $\pi_{i,j}(0)=p_iq_j$ for $1\leq i,j\leq N,$ and define
\begin{equation}\label{cconst0}
c_0 = \sum_{i,j=1}^N p_iq_iW_{i,j}
\end{equation}
If the marginal probabilities are different from zero, then since $\pi(0)$ is an interior point of $cP$ then $\by(0)\in A(\cP)$ and we can start our iterative procedure from there.
\subsection{Step decreasing procedure}
The sequential application of the maximum entropy method goes as follows.  The first step is to solve (\ref{prob2.1}) for the initial cost $c_0.$ Denote de maxentropic solution by $\bx_0.$ Then decrease $c$ and invoke MEM again. Whenever $\by(c_n)$ is in the relative interior of the data set (see Section \ref{remarks} in the appendix for more on this), then there is a $\bx_n^*$ satisfying $\bA\bx_n^*= \by(c_n).$ Suppose, for example, that the costs are decreased according to $c_n=c_{n-1}-\delta,$ for $n\geq 1$ with $c_0$ being the initial cost, and $\delta$ a small positive number. If we continue this process we shall arrive at a $n^*$ for which the corresponding $\by(c_{n*})$ is in the interior of $\bA(\cP),$ but $c_{n^*+1}$ is such that $\by(c_{n^*+1})$ is not in $\bA(\cP)$ and the maxentropic solution to (\ref{prob2a}) does not exist anymore. 

The process stops at this $n^*.$ The resulting optimal plan by $x^*.$ It will be of the type (\ref{sol1}), and it should be called the $\delta-$minimal cost plan. Observe that the difference between the $\delta-$minimal and the true minimum is less that $\delta.$

 \section{Appendices}

\subsection{Remarks about Problem \ref{prob2}}\label{remarks}

To make this self contained, here we cite some material developed in Gzyl (2020), put aside for not to interrupt the main discourse.
In the introduction we mention that it is convenient for notational purpose to relabel the unknowns in the problem. To a joint probability $\pi$ on $\{1,...,N\}^2$ we associate an $N^2-$vector $\bx\in \Omega.$ To express the marginality constraints in matrix form we consider an $(2N)\times N^2$ matrix $\bC$ constructed as follows:
\begin{definition}\label{constmat}
For $1\leq k\leq N$ the $k-$th row is the (transpose of) the $N^2-$vector with all component equal to $0$ except those at position $(k-1)N+j$ for $1\leq j\leq N$ which is equal to $1.$\\
For $N+1\leq k\leq 2N$ the $k-$th row is the (transpose of) the $N^2-$vector with zeros everywhere except at positions $(k-N)+(j-1)N$ for $n=(k-(N+1))N+j$ at which it equals $1$
\end{definition}
This association makes the dimension of the data space apparent. An example suffices to visualize the situation. For $X\in\{1,2,3\}$ and $Y\in\{1,2,3\}$ the constraint matrices look like:
\begin{example}\label{ex1}
\[ \bC = \left[ {\begin{array}{ccccccccc}
   1 & 1 & 1 & 0 & 0 & 0 & 0 &0 & 0 \\
   0 & 0 & 0 & 1 & 1 & 1 & 0 &0 & 0\\
   0 & 0 & 0 & 0 & 0 & 0 & 1 & 1 & 1 \\
   1 & 0 & 0  & 1 &  0 & 0 & 1 & 0 & 0 \\
   0 & 1 & 0 & 0 & 1 & 0 & 0 & 1 & 0\\
   0 & 0 & 1 & 0 & 0 & 1 & 0 & 0 & 1 \\
    \end{array} } \right].
\]
\end{example}

Since we require $\bx\in \Omega,$ the problem has convex constraints. Note that the sum of all row vectors of $\bC$ is $2\bu^t$ where $\bu$ is the $N^2-$vector of ones, and that the sum of the components of $\by$ is $2,$ then the constraint $\langle\bu,\bx\rangle=1$ is automatically satisfied if $\bx$ solves Problem (\ref{prob2}). That is the solution to (\ref{prob2}) actually lives in the simplex $\cP=\{\bx\in\Omega: \langle\bu,\bx\rangle=1\}.$ 

Let $\bC^{(1)}$ and $\bC^{(2)}$ denote the sub matrices of $\bC$ consisting, respectively, of the first $N$ and last $N$ rows. Note that:\\
{\bf (i)} As the rows of $\bC^{(1)}$ are independent vectors in $\bbR^{N^2},$ the rank of that sub-matrix is $N,$ therefore $\bC^{(1)}(\Omega)$ is a convex polytope in $\bbR^N$ with non empty interior. \\
{\bf (ii)} Similarly, the rows of  $\bC^{(2)}$  are independent vectors in $\bbR^{N^2}$ and the rank of that sub-matrix is $N,$ therefore $\bC^{(2)}(\Omega)$ is a convex polytope in $\bbR^N$ with non empty interior. 

Since the sum of the rows of $\bC^{(1)}$ and $\bC^{(2)}$ equals $\bu^t,$ the rows of the full matrix are not independent vectors. There are actually $2N-1$ independent rows among them. The rank of the constraint matrix is $2N-1.$ The image $\bC(\Omega)$ of $\Omega$ by $\bC$ is some convex $(2N-1)-$dimensional polytope in $\bbR^{2N},$ and the image $\cM=\bC(\cP)$ of $\cP$ is a polytope in $\bC(\Omega).$   The $(N^2-1)-$dimensional simplex $\cP$ is  the convex hull of the $N^2$ unit vectors $\be_n$ in $\bbR^{N^2}.$  To conclude\\
{\bf (iii)} Since $\bC$ is of rank $2N-1,$ then $\bC^{(1)}\bigcap\bC^{(2)}$ is a $1-$dimensional subspace of $\bbR^{2N}.$ Then $\bC(\Omega)$ is the convex sum  $\bC^{(1)}(\Omega)+\bC^{(2)}(\Omega),$ as well as 
$\bC(\cP)$ is the convex sum $\bC^{(1)}(\cP)+\bC^{(2)}(\cP),$ and therefore, the relative interior of $\bC(\cP)$ is not empty. When we add the positive row vector $\bW^t$ to $bC$ to obtain the matrix $\bA,$ we obtain $\bA(\cP)=\bC(\cP)\times[0,L]$ with $L=max\{W_n:n=1,...,N^2\}.$

\subsection{Appendix 2: Mathematical details about the maximum entropy method}\label{MD}
Here we describe in some detail the method of maximum entropy in the mean. Note that to solve Problem \ref{prob2a.1}, we need a solution to an ill-posed algebraic problem that satisfies the convexity constraint $\bxi:\Omega\to\Omega.$ The essence of MEM is  to think of the identity mapping $\bxi:\Omega\to\Omega$ as a random variable with respect to an unknown distribution. The maximum entropy part of the method comes in when determining a probability $\bbP$ upon $(\Omega,\cF)$ such that $\bA E_\bbP[\bxi]=\by.$ That is, the solution to the algebraic problem is the mean of $\bxi$ with respect to the probability that satisfies an entropy under some constraints.

To simplify the quest, it is convenient to start with some reference measure $\bbQ,$ upon $(\Omega,\cF)$ such that the convex hull $conv(supp(\bbQ))$ of its support equals the constraint set $\Omega,$ which is a closed, convex set. The other implicit aspect of our choice is that it makes computations as simple as possible. Since in our setup $\Omega$ is a product of closed intervals, a very convenient choice is:
\begin{equation}\label{meas}
\bbQ(d\xi) = \prod_{n=1}^{N^2}\Big(\varepsilon_{0}(d\xi_n)+\varepsilon_{1}(d\xi_n)\Big).
\end{equation}
We use the notation $\varepsilon_a$ to denote the measure that assigns unit point mass to point $a.$ The probability $\bbP\sim\bbQ$ is of the form
$$\bbP(d\xi) = \prod_{n=1}^{N^2}\Big(p_n\varepsilon_{0}(d\xi_n)+(1-p_n)\varepsilon_{1}(d\xi_n)\Big).$$
In general we would write $\bbP(d\xi) =\rho(\bxi)\bbQ(d\xi),$ which due to the special form of $\bbQ$ becomes the identity displayed above. Consider now the class
$$\cK = \{\bbP \sim \bbQ|\,\bA E_\bbP[\bxi] = \by(c)\}.$$
Any way to select a point from $\cK$ is valid a priori. A way that has proven to be very useful in many applications is by maximizing a very specific concave function defined on the class of all probabilities $\bbP\sim\bbQ$ by
\begin{definition}\label{ent}
$$S_\bbQ(\bbP) = -\int_\Omega\rho(\bxi)\ln\rho(\bxi)\bbQ(d\bxi) = -\sum_{n=1}^{N^2}p_n\ln p_n + (1-p_n)\ln(1-p_n),$$
called the entropy of $\bbP$ with respect to $\bbQ.$
\end{definition}
To explain the procedure to determine the $\blambda^*$ mentioned in Section 2, we will the following result proved by Kullback.
\begin{theorem}\label{kull}
Let $\bbP_1<<\bbQ$ and $\bbP_2<<\bbQ$ be two probabilities on $(\Omega,\cF),$  with densities $\rho_1(\bxi)$ and $\rho_2(\bxi).$ Define the (Kullback) divergence between $\bbP_1$ and $\bbP_2$ by
$$K(\rho_1,\rho_2) = \int_\Omega\rho_1(\bxi)\ln\big(\frac{\rho_1(\bxi)}{\rho_2(\bxi)}\big)\bbQ(d\bxi).$$
Then, $K(\rho_1,\rho_2)\geq 0$ and equals $0$ if and only if $\rho_1(\bxi)=\rho_2(\bxi)$ almost surely with respect to $\bbQ.$ 
\end{theorem}
We use this result as follows: Let $\bbP_1=\bbP \in\cK,$ and let $\bbP_2$ have exponential density given by
\begin{equation}\label{rep0}
\rho_{\blambda}(\bxi) = \frac{e^{-\langle\blambda,\bA\bxi\rangle}}{Z(\blambda)}
\end{equation}
\noindent where the normalization factor was defined in (\ref{norm1}) to be
$$Z(\lambda) = \int_\Omega e^{-\langle\blambda,\bA\bxi\rangle}d\bbQ(\bxi).$$
Invoking Theorem \ref{kull} we obtain
\begin{equation}\label{ineq1}
S_\bbQ(\bbP) \leq \ln Z(\blambda) + \langle\blambda,\by(c)\rangle \equiv \Sigma(\blambda,\by(c)).
\end{equation}
Therefore, the supremum on the left hand side is less or equal that the infimum of the right hand side. Thus if we could find a $\blambda^*$ such that  $\bbP^*=\rho_{\blambda^*}\bbQ \in \cK$we would have solved with the entropy maximization problem. This is what the standard method of maximum entropy is about. The result we want is
\begin{theorem}\label{main}
Let $\by(c)\in\bbR^{2N+1}$ be such that  $\Sigma(\blambda,\by(c))$ is bounded below. Then there is a unique $\blambda^*\in\bbR^{2N+1}$ such that $\bbP^*=\rho_{\blambda^*}\bbQ $ solves the entropy maximization problem, and
$$S_\bbQ(\bbP^*) = \ln Z(\blambda^*) + \langle\blambda^*,\by(c)\rangle = \Sigma(\blambda^*,\by(c)).$$
\end{theorem}
\begin{proof}
We have already mentioned all the necessary details. Note from (\ref{norm1}) that as $\{\blambda|Z(\blambda)<\infty\}=\bbR^{2N}.$ As the function $Z(\blambda)$ is continuously differentiable, the first order condition for $\blambda^*$ to be a minimum reads
$$\bA\int_\Omega \bxi\frac{e^{-\langle\blambda^*,\bA\bxi\rangle}}{Z(\blambda^*)}d\bbQ(\bxi) = \by(c)$$
From this we read off that $\bx^*$is given by given by (\ref{sol1}), and that is in the interior of the constraint set and their image is in the relative interior of $\bC(\Omega).$
\end{proof} 

\subsection{The solution to the cost minimization problem}
To tie the remarks made above to the cost minimization problem to the sequence of decreasing cost constrained problems, note that if $\by(c)$ is in the (relative) interior of $\bA\Omega,$ the optimal $\blambda^*$ exists. Otherwise, all that we can assert goes as follows.

For short, let us denote by $\Re$ the (relative) interior of the range of $\Omega$ by $\bA.$ Suppose that $\by(c)\in\partial\Re,$ and let $\by^{(k)}\in\Re\,\to \by(c)$ as $k \to \infty.$ Then there exist a sequence $\blambda^*_k\to \infty$ and $\bx^{(k)}$  as in Theorem \ref{main} such that $\bA\bx^{(k)}=\by^{(k)}$ and $\bx^{(k)}\to\bx^\infty\in\partial\Omega.$ Also
\begin{equation}\label{main2}
\bx^\infty_n = \left\{\begin{array}{c}
			a_n\;\;\;\mbox{whenever}\sum_{j=1}^{2N}\lambda^{(k)}_jA_{n,j}\to \infty\;\;\mbox{as}\;\;k\to\infty.\\
			b_n\;\;\;\mbox{whenever}\sum_{j=1}^{2N}\lambda^{(k)}_jA_{n,j}\to -\infty\;\;\mbox{as}\;\;k\to\infty.\\\end{array}\right.
\end{equation}
Not only that, as in our set up $\Omega$ is a finite dimensional (hyper-)box, we have as well that $\bA\bx^\infty=\by(c).$
The former comments can be formally stated as:
\begin{theorem}\label{main2}
If we know that the data vector $\by(c)\in\partial\Re,$ then there exists $\bx^\infty\in\partial\Omega$ such that $\bA\bx^\infty=\by(c).$
\end{theorem}
Obviously, if the data vector $\by\notin\bA(\Omega),$ no solution to the problem exists. In Section \ref{remarks} we saw that the problem of characterizing $\bA(\Omega)$ is not that simple, let alone characterizing its relative interior.

\subsection{Pending details for the geometry in pixel space}
Consider the mapping given in (\ref{corresp})between $int(\Omega)=\prod_{n=1}^{N^2}(a_n,b_n)$ and $\bbR^{N^2}$ given explicitly by
\begin{equation}\label{CoV}
\xi(\tau) = \frac{ae^{-a\tau}+be^{-b\tau}}{e^{-a\tau}+e^{-b\tau}}\;\;\Longleftrightarrow\;\;\tau=\frac{1}{D}\ln\bigg(\frac{b-\xi}{\xi-a}\bigg).
\end{equation} 
Due to the separability built into our problem we suppress reference to the label of the coordinates. Similarly, since the Hessian of $M(\btau)$ is diagonal and separable, when transporting the metric from (the tangent space to) $\bbR^{N^2}$ to (the tangent space of) $int(\Omega),$ we consider one generic coordinate. The pullback of  $M''$ back to is, by definition
\begin{equation}\label{PB}
g(\xi) = M''\big(\tau(\bx)\big)\bigg(\frac{d\tau}{d\xi}\bigg)^2\big.
\end{equation}
This definition ensures that the distance between two points does not depend on the system of coordinates used to describe the points. Put $D=b-a$. It is easy to see from (\ref{CoV}) that (\ref{PB}) becomes
\begin{equation}\label{PB1}
g(\xi) = \frac{1}{(b-\xi)(\xi-a)}.
\end{equation}
To find the geodesic distance between points in pixel space, we proceed as in Section 3 and put $g(\xi)=\big(k'(\xi)\big)^2,$ determine the function $k(\xi)$ to find the geodesic distance between points in pixel space. We shall follow two different routes.

{\bf First approach}\\
Put $S=(a+b)/2,$ add and subtract $S$ in each factor in the denominator of (\ref{PB1}) to obtain
$$g(\xi) =  \frac{1}{(\frac{D}{2})^2-(\xi-S)^2} = k'(\xi)^2,$$
that is we regard $k'(\xi)=\big[(\frac{D}{2})^2- (\xi-S)^2\big]^{-1/2}$ and therefore, since $a<\xi<b,$ we put:
$$k(\xi) = \int_a^\xi \frac{1}{\big[(\frac{D}{2})^2- (u-S)^2\big]^{1/2}} du.$$
Now do the following changes of variables: First put $v=\frac{2}{D}(u-S)$ and then $v=sin(\alpha)$  and $\theta(\xi)=\arcsin\big(\frac{2}{D}(\xi-S)\big)$ to obtain
\begin{equation}\label{CoV1}
k(\xi) = \int_{-1}^{\frac{2}{D}(u-S)}\frac{dv}{1-v^2} = \int_{-\pi/2}^{\theta(\xi)}d\alpha = 
\theta(\xi)+\frac{\pi}{2}.
\end{equation}
Now, put the coordinate labels back in place, let $G(\xi)$ be the matrix with elements $G_{n,m}(\bxi)=g_n(\xi_n)\delta_{n,m}.$ Then the definition is (\ref{PB} is such that for any continuously differentiable curve $\btau(t)$ in $\bbR^{N^2}$ its pullback $\bxi(t)=\xi(\btau(t))$ to $int(\Omega)=\prod_{n=1}^{N^2}(a_n,b_n)$ satisfies
\begin{equation}\label{basic1}
\langle\dot{\bxi},G(\xi)\dot{\bxi}\rangle = \langle\dot{\bxi},\big(\frac{\partial\btau}{\partial\bxi}\big)M''(\tau(\bxi))\big(\frac{\partial\btau}{\partial\bxi}\big)\dot{\bxi}\rangle = \langle\dot{\btau},M''(\btau)\dot{\btau}\rangle.
\end{equation}
We put $\big(\partial\btau/\partial\bxi\big)$ to denote the (symmetric) diagonal matrix with elements $d\tau_n/d\xi_n.$ Clearly $\dot{\btau}=\big(\partial\btau/\partial\bxi\big)\dot{\bxi}.$  What (\ref{basic1}) asserts is that the two curves have the same length. As a matter of fact, we have:
\begin{theorem}\label{eqdist}
Let $\tau(t)$ be a geodesic in the metric $M''$ such that $\btau(1)$ at $t=0$ and $\btau(2)$ at $t=1.$  Let $\bxi(t)=\xi(\btau(t))$ be its pullback with $\bxi(1)=\bxi(\btau(1))$ at $t=0$ and $\bxi(2)=\bxi(\btau(2))$ at $t=1.$ Then $\bxi(t)$ is a geodesic and 
\begin{eqnarray}\label{eqdist1}
d_G^2(\bxi(1),\bxi(2)) = \sum_{n=1}^{N^2}\big(k_n(\xi_n(1))-k_n(\xi_n(2)\big)^2\\
= \sum_{n=1}^{N^2}\big(h_n(\xi_n(1))-h_n(\xi_n(2)\big)^2 = d_M^2\big(\btau(1),\btau(2))\nonumber.
\end{eqnarray}
\end{theorem}
The proof hinges on the obvious remark that what (\ref{basic1}) means is that the length of the two curves is equal, that is:
$$\int_0^1 \bigg(\langle\dot{\bxi},G(\xi)\dot{\bxi}\rangle\bigg)^{1/2}dt = \int_0^1\bigg(\langle\dot{\btau},M''(\btau)\dot{\btau}\rangle\bigg)^{1/2}dt$$
Therefore, if the right hand side is minimal, so is the left hand side. The representation of the distances is as in Section 3.

{\bf Second approach}\\
As above,we use no coordinate labels. Consider the following computations
$$
k'(\xi) = \frac{1}{\big((\xi-a)(b-\xi)\big)^{1/2}}= \frac{\big((\xi-a)(b-\xi)\big)^{1/2}}{b-a}\big[\frac{1}{b-\xi} +\frac{1}{\xi-a}\big]
$$
We can write
$$\big((\xi-a)(b-\xi)\big)^{1/2}\frac{1}{b-\xi}=-2(\xi-a)\frac{1}{(\xi-a)^{1/2}}\frac{d(b-\xi)^{1/2}}{d\xi}.$$
Similarly
$$\big((\xi-a)(b-\xi)\big)^{1/2}\frac{1}{\xi-a} = -2(\xi-a)(b-\xi)^{1/2}\frac{d(\xi-a)^{-1/2}}{d\xi}.$$
Putting this together we have
$$
k(\xi) = \int_a^\xi\frac{du}{\big((u-a)(b-u)\big)^{1/2}} = 2\int_a^\xi\frac{(u-a)}{1+\frac{(b-u)}{(u-a)}}\frac{d}{du}\bigg(\frac{(b-u)}{(u-a)}\bigg)^{1/2}.$$
That is,
$$k(\xi) = \int_a^\xi \frac{d}{du}\arctan\bigg[\bigg(\frac{b-u}{u-a}\bigg)^{1/2}\bigg]du = 2\arctan\bigg[\bigg(\frac{b-\xi}{\xi-a}\bigg)^{1/2}\bigg] + \frac{\pi}{2}.$$

To conclude, form (\ref{CoV}) we obtain $\exp(D\tau/2) = \big((b-\xi)/(\xi-a)\big)^{1/2},$ and therefore $k(\xi) = h(\tau(\xi)).$ This is a longer way of saying that the factorizations of $G$ and $M''$ are consistent in the sense that $k'(\xi)=h'(\tau(\xi))d\tau/d\xi.$ The roundabout way of getting at $k(\xi)=h(\tau(\xi)$ was to relate $k(\xi)$ to the factorization of the pullback metric $G$ in order to establish the identification of $d_M(\btau(1),\btau(2))=d_G(\bxi(1),\bxi(2)$ as geodesic distances.

\end{document}